\documentclass[12pt,oneside]{article}
\usepackage{amsfonts,pstricks}
\usepackage{amssymb}
\usepackage{euscript}
\usepackage{eufrak}
\usepackage{mathrsfs}
\usepackage{latexsym}
\usepackage{amsmath,enumerate}

\newtheorem{thm}{Theorem}[section]
\newtheorem{lem}{Lemma}[section]

\newtheorem{rmk}{Remark}[section]
\newtheorem{defn}{Definition}[section]
\newtheorem{prop}{Proposition}[section]
\newtheorem{ex}{Example}[section]

\newcommand{\NI}{\noindent}
\newcommand{\eb}{\vrule width 0.22 true cm height 0.22 true cm depth 0pt}

\usepackage{graphicx}
\date{}
\usepackage[normalem]{ulem}
\usepackage{hyperref}

\begin{document}
	
	\title
	{\large \bf Characterizations of $e^\star$-open sets and nearby open sets  on Infra topological spaces}
	\author{G. Saravanakumar$^1$\footnote{saravananguru2612@gmail.com}\hspace{.2cm} and D. Jaya Priya$^2$\footnote{jayapriyadhanasekar001@gmail.com}   \footnote{$^{1,2}$Department of Mathematics, Vel Tech Rangarajan Dr.Sagunthala R\&D Institute of Science and Technology (Deemed to be University), Avadi, Chennai-600062, India. }}
	\maketitle

	\begin{abstract}
		\NI \indent The study of infra-topological spaces focuses on characterizations of $e^\star$-open sets and nearby open sets in infra-topological spaces. The $e^\star$-open sets, a variation of open sets, are explored for their unique properties and relationships within the infra-topological framework. Additionally, nearby open sets, which capture the notion of points being close to each other, are investigated to provide a comprehensive understanding of the topological structure. The research aims to contribute to the broader field of topology by extending traditional concepts to infra-topological spaces, offering new perspectives on openness and proximity. The findings not only deepen our understanding of mathematical structures but also open avenues for applications in various scientific and engineering disciplines.
	\end{abstract}
	
	
	\textbf{AMS (2000) subject classification:}54A05,54A10

\section{Introduction}

A seminal work by S. Mashhour et al. \cite{mas} laid the foundation for supra topological spaces, wherein they delved into the intricacies of $s$-continuous functions and $s^{*}$-continuous functions. The groundbreaking exploration of supra topological spaces marked a significant contribution to the field.

Building upon this foundation, Adel. M. Al. Odhari \cite{ade} extended the theoretical landscape by introducing and studying infra topological spaces. Odhari's work included a comprehensive examination of open sets within infra topological spaces, shedding light on their fundamental properties.

E. Ekici, in a series of noteworthy contributions \cite{eki1,eki2,eki3,eki4,eki5}, introduced and extensively investigated the properties of $e$ and $e^\star$ along with nearby open sets in the context of general topological spaces. Ekici's research provided valuable insights into the behavior of these sets, enriching the understanding of topological structures.

Motivated by these advancements, our present work focuses on extending the study of $e$ and $e^\star$ open sets to infra topological structures. We embark on an in-depth investigation, exploring the properties and characteristics of these sets within the specific framework of infra topological spaces. To bolster our findings, we present concrete examples that not only validate our assumptions but also serve as illustrative instances of the nuanced interplay between infra topological structures and the introduced open sets. Through this research, we aim to contribute to the ongoing dialogue surrounding the interplay of different topological spaces and the behavior of distinct open sets within them.

\section{Preliminaries}

\defn\em\label{2.1} Let $X$ be any arbitrary set. An infra topological space ($\mathcal{ITS}$) on $X$ is a collection infra topology ($\mathcal{IT}$) of subsets of $X$ such that the following axioms are satisfying:

\begin{enumerate}[(i)]
	\item $\emptyset,~X\in\mathcal{IT},.$
	\item The intersection of the elements of any subcollection of $\mathcal{IT}$ in $X.$
\end{enumerate}
i.e, If $O_{i}\in\mathcal{IT},~1\subseteq i\subseteq n,~\cap O_{i}\in\mathcal{IT}.$

Terminology, the ordered pair $(X,\mathcal{IT})$ is called  $\mathcal{ITS}$. we simply say $X$ is a $\mathcal{ITS}$. 

\defn\em Let $(X,~\mathcal{IT})$ be an $\mathcal{ITS}$ and $A\subset X.~~ A$ is called  infra open set ($\mathcal{IOS}$) if $A\in\mathcal{IT}.$ The complement of $\mathcal{IOS}$ is called infra closed set ($\mathcal{ICS}$). 

\defn\em Let $(X,~\mathcal{IT})$ be $\mathcal{ITS}$. A subset $C\subset X$ is called $\mathcal{ICS}$ in $X$ if $X-C$ is  $\mathcal{IOS}$ in $X$. That is $C$ is  $\mathcal{ICS}$ iff $X-C\in\mathcal{IT}.$



\thm\em Let $(X,~\mathcal{IT})$ be a Topological Space($\mathcal{TS}$) , then $(X,~\mathcal{IT})$ is an  $\mathcal{ITS}$.




\thm\em Let $(X,~\mathcal{IT})$ be  $\mathcal{ITS}$. Then:
\begin{enumerate}[(i)]
	\item $\emptyset,~X$ are  $\mathcal{IOS}$.
	\item Any arbitrary intersections of $\mathcal{IOS}'$s are $\mathcal{IOS}'$s.
\end{enumerate}

%
%
%

\thm\em let $(X,~\mathcal{IT})$ and $(X,~\mathcal{IT}^{*})$ be two  $\mathcal{ITS}'$s on set $X$. Then the intersection $\mathcal{IT}$ and $\mathcal{IT}^{*}$ is an $\mathcal{ITS}$, while the union $\mathcal{IT}$ and $\mathcal{IT}^{*}$ not necessarily.

\defn\em Let $(X,~\mathcal{IT})$ be an $\mathcal{ITS}$ and $A\subset X$. A point $x\in X$ is called Infra-Cluster Point ($\mathcal{IC_P}$) of $A$, if for all $\mathcal{IOS}$ $O$ containing $x$, then $A\cap(O\backslash \{x\})\neq\emptyset.$

\defn\em Let $(X,~\mathcal{IT})$ be an $\mathcal{ITS}$ and $A\subset X$. The set of all  $\mathcal{IC_P}$ of $A$ is called the  Infra Derived Set ($\mathcal{IDS}$) of $A.$

\thm\em Let $(X,~\mathcal{IT})$ be  $\mathcal{ITS}$. Then:
\begin{enumerate}[(i)]
	\item $\emptyset,~X\in\mathcal{IT}$ are $\mathcal{ICS}$.
	\item Any arbitrary finite intersections of $\mathcal{ICS}'$s is an $\mathcal{ICS}'$s.
\end{enumerate}


\defn\em Let $(X,~\mathcal{IT})$ be an $\mathcal{ITS}$ and $A\subset X$.  The infra closure  ($\mathcal{I.CL}$) of $A$ is a set denoted by $\mathcal{I.CL}(A)$ and given by: $\mathcal{I.CL}(A)=\bigcap\{C_{i}:A\subset C_{i},~X-C_{i}\in\mathcal{IT}\}.$That is, $\mathcal{I.CL}(A)$ is the intersection of all $\mathcal{ICS}$ contained the set $A.$

\rmk\em Since $\mathcal{I.CL}(A)$ is the intersection of all $\mathcal{ICS}'$s containing in $A$, then $A\subset \mathcal{I.CL}(A)$ and $\mathcal{I.CL}(A)$ is the smallest $\mathcal{ICS}'$s.

\defn\em Let $(X,~\mathcal{IT})$ be an $\mathcal{ITS}$ and $A\subset X$.  The Infra Interior ($\mathcal{I.INT}$) of $A$ is a set denoted by $\mathcal{I.INT}(A)$ and given by: $\mathcal{I.INT}(A)=\bigcup\{O_{i}:O_{i}\subset A,~O_{i}\in\mathcal{IT}\}.$That is, $\mathcal{I.INT}(A)$ is the union of all $\mathcal{IOS}$ contained in the set $A.$

\rmk\em Since $\mathcal{I.INT}(A)$ is the union of all $\mathcal{IOS}$'s contained in $A$, then $\mathcal{I.INT}(A)\subset A$ and $\mathcal{I.CL}(A)$ is the smallest $\mathcal{IOS}'$s. Also if $O$ is $\mathcal{IOS}$ contained in $A$, then $O\subset \mathcal{I.INT}(A)$ .

\defn\em Let $(X,~\mathcal{IT})$ be an $\mathcal{ITS}$ and $A\subset X$.  The Infra Exterior Point ($\mathcal{IEP}$) of $A$ is a set denoted by $\mathcal{IEP}(A)$ and given by: $\mathcal{IEP}(A)=\mathcal{I.INT}(A^{c}).$That is, Set of all $\mathcal{I.INT}$ of complement of $A.$

\defn\em Let $(X,~\mathcal{IT})$ be an $\mathcal{ITS}$ and $A\subset X$.  The Infra-Boundary Points ($\mathcal{IBP}$) of $A$ is a set denoted by $\mathcal{IBP}(A)$ and given by: $\mathcal{IBP}(A)=X\backslash \mathcal{I.INT}(A)\cup \mathcal{IEP}(A)$

\thm\em Let $(X,~\mathcal{IT})$ be an $\mathcal{ITS}$ and $A,~B\subset X.$ The $\mathcal{IDS}$ Axioms satisfies the followings:
\begin{enumerate}[(i)]
	\item  $\mathcal{IDS}(\emptyset)=\emptyset.$
	\item  If $A\subset B$ then $\mathcal{IDS}(A)\subset \mathcal{IDS}(B)$
	\item If $x\in \mathcal{IDS}(A)$ then $x\in \mathcal{IDS}(A\backslash \{x\})$ . 
	\item $\mathcal{IDS}(A\cap B)\subset \mathcal{IDS}(A)\cap \mathcal{IDS}(B)$.
	\item $\mathcal{IDS}(A\cup B)=\mathcal{IDS}(A)\cup \mathcal{IDS}(B)$ .
\end{enumerate}

\thm\em Let $(X,~\mathcal{IT})$ be an $\mathcal{ITS}$ and $A,B\subset X$.The  $\mathcal{I.CL}$ Axioms satisfying the following conditions:
\begin{enumerate}[(i)]
	\item $A$ is $\mathcal{ICS}$ iff  $A=\mathcal{I.CL}(A)$ .
	\item $ \mathcal{I.CL}(\emptyset)=\emptyset$  and $\mathcal{I.CL}(X)=X.$
	\item $\mathcal{I.CL}(\mathcal{I.CL}(A))=\mathcal{I.CL}(A)$ .
	\item If $A\subset B$ then $\mathcal{I.CL}(A)\subset \mathcal{I.CL}(B)$ .
	\item $\mathcal{I.CL}(A\cap B)\subset \mathcal{I.CL}(A)\cap \mathcal{I.CL}(B)$ .
\end{enumerate}

\thm\em Let $(X,~\mathcal{IT})$ be an $\mathcal{ITS}$ and $A,B\subset X.$ The $\mathcal{I.INT}$ Axioms  given by:
\begin{enumerate}[(i)]
	\item $A$ is $\mathcal{IOS}$ iff  $A=\mathcal{I.INT}(A)$ . 
	\item $ \mathcal{I.INT}(X)=X$ and $\mathcal{I.INT}(\emptyset)=\emptyset.$
	\item $\mathcal{I.INT}(\mathcal{I.INT}\subset(A))=\mathcal{I.INT}(A)$ .
	\item If $A\subset B$ then $\mathcal{I.INT}(A)\subset \mathcal{I.INT}(B)$ . 
	\item $\mathcal{I.INT}(A\cap B)=\mathcal{I.INT}(A)\cap \mathcal{I.INT}(B)$ .
\end{enumerate}

\thm\em Let $(X\mathcal{IT})$ be an $\mathcal{ITS}$ and $A,B\subset X.$ The   $\mathcal{IEP}$ Axioms  given by:
\begin{enumerate}[(i)]
	\item $\mathcal{IEP}(X)=\emptyset$ and $\mathcal{IEP}(\emptyset)=X.$
	\item $\mathcal{IEP}(A)\subset A^{c}.$
	\item $\mathcal{IEP}(A\cup B)=\mathcal{IEP}(A)\cap \mathcal{IEP}(B)$ . 
	\item If $A\subset B,$ then $\mathcal{IEP}(B)\subset \mathcal{IEP}(A)$ . 
	\item $\mathcal{IEP}(A\cap B)\subset \mathcal{IEP}(A)\cup \mathcal{IEP}(B)$ .
\end{enumerate}

\thm\em Let $(X,~\mathcal{IT})$ be an $\mathcal{ITS}$ and $A\subset X.$ The  $\mathcal{IB}$ Axioms given by:
\begin{enumerate}[(i)]
	\item $\mathcal{IBP}(X)=\mathcal{IBP}(\emptyset)=\emptyset.$
	\item $\mathcal{IBP}(A\cap B)=\mathcal{IBP}(A)\cup \mathcal{IBP}(B)$ .
\end{enumerate}

\thm\em Let $(X,~\mathcal{IT})$ be an (ITS) and $A\subset X.$ Then: 
\begin{enumerate}[(i)]
	\item $A\subset \mathcal{I.CL}(A)\rightarrow \mathcal{IDS}(A)\subset \mathcal{IDS}(\mathcal{I.CL}(A))$.
	\item $\mathcal{I.INT}(A)\subset A\rightarrow \mathcal{IDS}(\mathcal{I.INT}(A))\subset \mathcal{IDS}(A)$ .
	\item If $A$ is $\mathcal{ICS}$, then $\mathcal{IDS}(A)\subset A.$
	\item $\mathcal{I.CL}(A)=A\cup \mathcal{IDS}(A)$ .
	\item  $\mathcal{IB}(A)=\mathcal{I.CL}(A)\backslash \mathcal{I.INT}(A)$ .
	\item  $\mathcal{I.CL}(A)=\mathcal{IBP}(A)\cup \mathcal{I.INT}(A)$ .
	\item  $\mathcal{IBP}(A)\subset \mathcal{I.CL}(A)$ .
	\item  $\mathcal{I.INT}(A)\cap \mathcal{IBP}(A)=\emptyset.$
\end{enumerate}

\section{Characterizations of $e$-open and nearby open sets sets}

\defn\em Let $(X,~\mathcal{IT})$ be an $\mathcal{ITS}$ and $A\subset X$.  The Infra $\delta$ Interior ($\mathcal{I.INT_\delta}$) of $A$ is a set denoted by $\mathcal{I.INT_\delta}(A)$ and given by: $\mathcal{I.INT_\delta}(A)=\bigcup\{O_{i}:O_{i}\subset A,~O_{i}\in\mathcal{I.}\mathcal{ROS}\}.$

\defn\em Let $(X,~\mathcal{IT})$ be an $\mathcal{ITS}$ and $A\subset X$.  The infra $\delta$ closure  ($\mathcal{I.CL_\delta}$) of $A$ is a set denoted by $\mathcal{I.CL_\delta}(A)$ and given by: $\mathcal{I.CL_\delta}(A)=\bigcap\{C_{i}:A\subset C_{i},~X-C_{i}\in\mathcal{I.}\mathcal{RCS}\}.$

\defn\em Let $(X,~\mathcal{IT})$ be  $\mathcal{ITS}$. Then a Set $A$ is said to be Infra 
\begin{enumerate}[(i)]
	\item regular open set($\mathcal{I.}\mathcal{ROS}$) if $A=\mathcal{I.INT}(\mathcal{I.CL}(A))$
	\item $\delta$-open set ($\mathcal{I.}\delta\mathcal{OS}$) if $A=\mathcal{I.INT}_\delta(A)$
	\item pre open set($\mathcal{IPOS}$) if $A\subseteq\mathcal{I.INT}(\mathcal{I.CL}(A))$
	\item semi open set($\mathcal{ISOS}$) if $A\subseteq\mathcal{I.CL}(\mathcal{I.INT}(A))$
	\item $\delta$-pre open set($\mathcal{I.}\delta\mathcal{POS}$) if $A\subseteq\mathcal{I.INT}(\mathcal{I.CL}_\delta(A))$
	\item $\delta$-semi open set($\mathcal{I.}\delta\mathcal{SOS}$) if $A\subseteq\mathcal{I.CL}(\mathcal{I.INT}_\delta(A))$
	\item $e$-open set($\mathcal{I.}e\mathcal{OS}$) if $A\subseteq \mathcal{I.CL}(\mathcal{I.INT}_\delta(A))\cup \mathcal{I.INT}(\mathcal{I.CL}_\delta(A))$
	\item $e^*$-open set ($\mathcal{I.}e^*\mathcal{OS}$) if $A\subseteq \mathcal{I.CL}(\mathcal{I.INT}(\mathcal{I.CL}_\delta (A)))$
	\item $a$-open set ($\mathcal{I.}a\mathcal{OS}$) if $A\subseteq \mathcal{I.INT}(\mathcal{I.CL}(\mathcal{I.INT}_\delta (A)))$
	\item $\beta$-open set($\mathcal{I.}\beta\mathcal{OS}$) if $A\subseteq \mathcal{I.CL}(\mathcal{I.INT}(\mathcal{I.CL}(A)))$
\end{enumerate}

\defn\em Let $(X,~\mathcal{IT})$ be $\mathcal{ITS}$. A subset $C\subset X$ is called  $\mathcal{I.}\mathcal{RCS},(\text{resp.}\mathcal{I.}\delta\mathcal{CS},\mathcal{IPCS},\mathcal{ISCS},\mathcal{I.}\delta\mathcal{PCS},\mathcal{I.}\delta\mathcal{SCS},\mathcal{I.}a\mathcal{CS},\mathcal{I.}e\mathcal{CS},\mathcal{I.}e^*\mathcal{CS},\mathcal{I.}\beta\mathcal{CS}$) in $X$ if $X-C$ is $\mathcal{I.}\mathcal{ROS},(\text{resp.}\mathcal{I.}\delta\mathcal{OS},\mathcal{IPOS},\mathcal{ISOS},\mathcal{I.}\delta\mathcal{POS},\mathcal{I.}\delta\mathcal{SOS},$ $\mathcal{I.}a\mathcal{OS},\mathcal{I.}e\mathcal{OS},$ $\mathcal{I.}e^*\mathcal{OS},\mathcal{I.}\beta\mathcal{OS}$) in $X$.

\defn\em Let $(X,~\mathcal{IT})$ be an  and $A\subset X$. A point $x\in X$ is called Infra

\begin{enumerate}[(i)]
	\item $\delta$ Cluster Point($\mathcal{I.}\delta\mathcal{C_P}$) if for all $\mathcal{I.}\delta\mathcal{OS}$ $O$ containing $x$, then $A\cap(O\backslash \{x\})\neq\emptyset.$
	\item Pre Cluster Point($\mathcal{IPC_P}$) if for all $\mathcal{IPOS}$ $O$ containing $x$, then $A\cap(O\backslash \{x\})\neq\emptyset.$
	\item  Semi Cluster Point($\mathcal{ISC_P}$) if for all $\mathcal{ISOS}$ $O$ containing $x$, then $A\cap(O\backslash \{x\})\neq\emptyset.$
	\item $\delta$Pre Cluster Point($\mathcal{I.}\delta\mathcal{PC_P}$) if for all $\mathcal{I.}\delta\mathcal{POS}$ $O$ containing $x$, then $A\cap(O\backslash \{x\})\neq\emptyset.$
	\item $\delta$Semi Cluster Point($\mathcal{I.}\delta\mathcal{SC_P}$) if for all $\mathcal{I.}\delta\mathcal{SOS} $ $O$ containing $x$, then $A\cap(O\backslash \{x\})\neq\emptyset.$
	\item $a$ Cluster Point($\mathcal{I.}a\mathcal{C_P}$) if for all $\mathcal{I.}a\mathcal{OS}$ $O$ containing $x$, then $A\cap(O\backslash \{x\})\neq\emptyset.$
	\item $e$ Cluster Point($\mathcal{I.}e\mathcal{C_P}$) if for all $\mathcal{I.}e\mathcal{OS}$ $O$ containing $x$, then $A\cap(O\backslash \{x\})\neq\emptyset.$
	\item $e^*$ Cluster Point($\mathcal{I.}e^*\mathcal{C_P}$) if for all $\mathcal{I.}e^*\mathcal{OS}$ $O$ containing $x$, then $A\cap(O\backslash \{x\})\neq\emptyset.$
	\item $\beta$ Cluster Point($\mathcal{I.}\beta\mathcal{C_P}$) if for all $\mathcal{I.}\beta\mathcal{OS}$ $O$ containing $x$, then $A\cap(O\backslash \{x\})\neq\emptyset.$
\end{enumerate}

\defn\em Let $(X,~\mathcal{IT})$ be an (ITS) and $A\subset X$. The set of all $\mathcal{I.}\delta\mathcal{C_P},(\text{resp.}\mathcal{IPC_P},\mathcal{ISC_P},\mathcal{I.}\delta\mathcal{PC_P},\mathcal{I.}\delta\mathcal{SC_P},\mathcal{I.}a\mathcal{C_P},\mathcal{I.}e\mathcal{C_P},\mathcal{I.}e^*\mathcal{C_P},\mathcal{I.}\beta\mathcal{C_P}$)  of $A$ is called the Infra $\delta$ Derived Set(resp. Infra Pre Derived Set,Infra Semi Derived Set,Infra $\delta$Pre Derived Set, Infra $\delta$Semi Derived Set,Infra $a$ Derived Set,Infra $e$ Derived Set,Infra $e^*$ Derived Set,Infra $\beta$ Derived Set)  of $A$  and is denoted by $\mathcal{I.}\delta \mathcal{DS}(resp. (\text{resp.}\mathcal{IPDS},\mathcal{ISDS},\mathcal{I.}\delta\mathcal{PDS},\mathcal{I.}\delta\mathcal{SDS},\mathcal{I.}a\mathcal{DS},\mathcal{I.}e\mathcal{DS},\mathcal{I.}e^*\mathcal{DS},\mathcal{I.}\beta\mathcal{DS}))$ of $A$.

\thm\em Let $(X,~\mathcal{IT})$ be  $\mathcal{ITS}$. Then:
\begin{enumerate}[i.]
	\item $\emptyset,~X\in\mathcal{IT}$ are  $\mathcal{ICS}$.
	\item Any arbitrary finite intersections of $\mathcal{ICS}'$s is an $\mathcal{ICS}'$s.
	\end{enumerate}
\textbf{Proof:}
(i) Since $X-\emptyset=X\in\mathcal{IT}$  and $X-X=\emptyset\in\mathcal{IT}$ are $\mathcal{ICS}'$s.

(ii) Let $\{C_{i}\ :\ i\in I\}$ be an arbitrary family of $\mathcal{ICS}'$s such that $C_{i}\in\mathcal{IT}$ for all $i\in I.$ Now, $X-C_{i}\in\mathcal{IT}$ is $\mathcal{IOS}$ for all $i\in I.$ But $X-C_{i}=C_{i}^{c}\in\mathcal{IT}$ then $\bigcap C_{i}^{c}=\bigcap(X-C_{i})=X-\bigcap C_{i}\in\mathcal{IT},\forall i\in I$. Hence $\bigcap C_{i}\in\mathcal{IT},\forall i\in I$ is $\mathcal{ICS}$.
 

\defn\em Let $(X,~\mathcal{IT})$ be an $\mathcal{ITS}$ and $A\subset X$.  The Infra 
\begin{enumerate}[(i)]
	\item Pre Closure ($\mathcal{I.PCL}$) of $A$ and given by $\mathcal{I.PCL}(A)=\bigcap\{C_{i}:A\subset C_{i},~X-C_{i}\in \mathcal{IPOS} \}.$
	\item Semi Closure ($\mathcal{I.SCL}$) of $A$ and given by $\mathcal{I.SCL}=\bigcap\{C_{i}:A\subset C_{i},~X-C_{i}\in \mathcal{ISOS} \}.$
	\item $\delta$Pre Closure ($\mathcal{I.}\delta\mathcal{PCL}$) of $A$ and given by $\mathcal{I.}\delta\mathcal{PCL}=\bigcap\{C_{i}:A\subset C_{i},~X-C_{i}\in \mathcal{I.}\delta\mathcal{POS} \}.$
	\item $\delta$Semi Closure ($\mathcal{I.}\delta\mathcal{SCL}$) of $A$ and given by $\mathcal{I.}\delta\mathcal{SCL}=\bigcap\{C_{i}:A\subset C_{i},~X-C_{i}\in \delta\mathcal{SOS} \}.$
	\item $a$ Closure ($\mathcal{I.}a\mathcal{CL}$) of $A$ and given by $\mathcal{I.}a\mathcal{CL}=\bigcap\{C_{i}:A\subset C_{i},~X-C_{i}\in \mathcal{I.}a\mathcal{OS} \}.$
	\item $e$ Closure ($\mathcal{I.}e\mathcal{CL}$) of $A$ and given by $\mathcal{I.}e\mathcal{CL}=\bigcap\{C_{i}:A\subset C_{i},~X-C_{i}\in \mathcal{I.}e\mathcal{OS} \}.$
	\item $e^*$ Closure ($\mathcal{I.}e^*\mathcal{CL}$) of $A$ and given by $\mathcal{I.}e^*\mathcal{CL}=\bigcap\{C_{i}:A\subset C_{i},~X-C_{i}\in \mathcal{I.}e^*\mathcal{OS} \}.$
	\item $\beta$ Closure ($\mathcal{I.}\beta\mathcal{CL}$) of $A$ and given by $\mathcal{I.}\beta\mathcal{CL})=\bigcap\{C_{i}:A\subset C_{i},~X-C_{i}\in \mathcal{I.}\beta\mathcal{OS}) \}.$
	\end{enumerate}


\defn\em Let $(X,~\mathcal{IT})$ be an $\mathcal{ITS}$ and $A\subset X$.  The Infra 
\begin{enumerate}[(i)]
	\item Pre Interior ($\mathcal{I.PINT}$) of $A$ and given by $\mathcal{I.PINT}(A)=\bigcup\{C_{i}:A\supset C_{i},~C_{i}\in \mathcal{IPOS} \}.$
	\item Semi Interior ($\mathcal{I.SINT}$) of $A$ and given by $\mathcal{I.SINT}=\bigcup\{C_{i}:A\supset C_{i},~C_{i}\in \mathcal{ISOS} \}.$
	\item $\delta$Pre Interior ($\mathcal{I.}\delta\mathcal{PNT}$) of $A$ and given by $\mathcal{I.}\delta\mathcal{PINT}=\bigcup\{C_{i}:A\supset C_{i},~C_{i}\in \mathcal{I.}\delta\mathcal{POS} \}.$
	\item $\delta$Semi Interior ($\mathcal{I.}\delta\mathcal{SINT}$) of $A$ and given by $\mathcal{I.}\delta\mathcal{SINT}=\bigcup\{C_{i}:A\supset C_{i},~C_{i}\in \delta\mathcal{SOS} \}.$
	\item $a$ Interior($\mathcal{I.}a\mathcal{INT}$) of $A$ and given by $\mathcal{I.}a\mathcal{INT}=\bigcup\{C_{i}:A\supset C_{i},~C_{i}\in \mathcal{I.}a\mathcal{OS} \}.$
	\item $e$ Interior ($\mathcal{I.}e\mathcal{INT}$) of $A$ and given by $\mathcal{I.}e\mathcal{INT}=\bigcup\{C_{i}:A\supset C_{i},~C_{i}\in \mathcal{I.}e\mathcal{OS} \}.$
	\item $e^*$ Interior ($\mathcal{I.}e^*\mathcal{INT}$) of $A$ and given by $\mathcal{I.}e^*\mathcal{INT}=\bigcup\{C_{i}:A\supset C_{i},~C_{i}\in \mathcal{I.}e^*\mathcal{OS} \}.$
	\item $\beta$ Interior ($\mathcal{I.}\beta\mathcal{INT}$) of $A$ and given by $\mathcal{I.}\beta\mathcal{INT})=\bigcup\{C_{i}:A\supset C_{i},~C_{i}\in \mathcal{I.}\beta\mathcal{OS}) \}.$
\end{enumerate}

\defn\em Let $(X,~\mathcal{IT})$ be an $\mathcal{ITS}$ and $A\subset X$.  The Infra $\delta$ Exterior (resp. Infra Pre Exterior,Infra Semi Exterior,Infra $\delta$Pre Exterior, Infra $\delta$Semi Exterior,Infra $a$ Exterior,Infra $e$ Exterior,Infra $e^*$ Exterior,Infra $\beta$ Exterior)  of $A$ is a set denoted by $\mathcal{I.}\delta\mathcal{EX}(\text(\text{resp.}~\mathcal{IPEX},$ $\mathcal{ISEX},\mathcal{I.}\delta\mathcal{PEX},\mathcal{I.}\delta\mathcal{SEX},\mathcal{I.}a\mathcal{EX},\mathcal{I.}e\mathcal{EX},\mathcal{I.}e^*\mathcal{EX},\mathcal{I.}\beta\mathcal{EX})$ of $A$ and given by: $\mathcal{I.}\delta\mathcal{EX}(\text{resp.}~\mathcal{IPEX},\mathcal{ISEX},\mathcal{I.}\delta\mathcal{PEX},\mathcal{I.}\delta\mathcal{SEX},\mathcal{I.}a\mathcal{EX},\mathcal{I.}e\mathcal{EX},\mathcal{I.}e^*\mathcal{EX},\mathcal{I.}\beta\mathcal{EX})(A)$ $=\mathcal{I.}\delta\mathcal{INT}(\text(\text{resp.}~\mathcal{IPIP},\mathcal{ISIP},\mathcal{I.}\delta\mathcal{PIP},\mathcal{I.}\delta\mathcal{SIP},\mathcal{I.}a\mathcal{INT},\mathcal{I.}e\mathcal{INT},\mathcal{I.}e^*\mathcal{INT},$ $\mathcal{I.}\beta\mathcal{INT})(A^c).$

\defn\em Let $(X,~\mathcal{IT})$ be an $\mathcal{ITS}$ and $A\subset X$. The Infra
\begin{enumerate}[(i)]
	\item   $\delta$ Boundary (briefly,$\mathcal{I.}\delta\mathcal{B})$ of $A$ is described and indicated by $\mathcal{I.}\delta\mathcal{B}(A)=X\backslash \mathcal{I.}\delta \mathcal{INT}(A)\cup \mathcal{I.}\delta \mathcal{EX}(A)$
	\item   Semi Boundary (briefly,$\mathcal{I.}\mathcal{SB})$ of $A$ is described and indicated by $\mathcal{I.}\mathcal{SB}(A)=X\backslash \mathcal{I.} \mathcal{SINT}(A)\cup \mathcal{I.}\mathcal{SEX}(A)$
	\item  pre Boundary (briefly,$\mathcal{I.}\mathcal{PB})$ of $A$ is described and indicated by $\mathcal{I.}\mathcal{PB}(A)=X\backslash \mathcal{I.} \mathcal{PINT}(A)\cup \mathcal{I.} \mathcal{PEX}(A)$
	\item  $\delta$ Semi Boundary (briefly,$\mathcal{I.}\delta\mathcal{SB})$ of $A$ is described and indicated by $\mathcal{I.}\delta\mathcal{SB}(A)=X\backslash \mathcal{I.}\delta \mathcal{SINT}(A)\cup \mathcal{I.}\delta \mathcal{SEX}(A)$
	\item  $\delta$ pre Boundary (briefly,$\mathcal{I.}\delta\mathcal{PB})$ of $A$ is described and indicated by $\mathcal{I.}\delta\mathcal{PB}(A)=X\backslash \mathcal{I.}\delta \mathcal{PINT}(A)\cup \mathcal{I.}\delta \mathcal{PEX}(A)$
	\item  $a$ Boundary (briefly,$\mathcal{I.}a\mathcal{B})$ of $A$ is described and indicated by $\mathcal{I.}a\mathcal{B}(A)=X\backslash \mathcal{I.}a \mathcal{INT}(A)\cup \mathcal{I.}a \mathcal{EX}(A)$
	\item  $e$ Boundary (briefly,$\mathcal{I.}e\mathcal{B})$ of $A$ is described and indicated by $\mathcal{I.}e\mathcal{B}(A)=X\backslash \mathcal{I.}e \mathcal{INT}(A)\cup \mathcal{I.}e \mathcal{EX}(A)$
	\item  $e^*$ Boundary (briefly,$\mathcal{I.}e^*\mathcal{B})$ of $A$ is described and indicated by $\mathcal{I.}e^*\mathcal{B}(A)=X\backslash \mathcal{I.}e^* \mathcal{INT}(A)\cup \mathcal{I.}e^* \mathcal{EX}(A)$
	\item  $\beta$ Boundary (briefly,$\mathcal{I.}\beta\mathcal{B})$ of $A$ is described and indicated by $\mathcal{I.}\beta\mathcal{B}(A)=X\backslash \mathcal{I.}\beta \mathcal{INT}(A)\cup \mathcal{I.}\beta \mathcal{EX}(A)$
\end{enumerate}

\thm\em\label{t3.2} Let $(X,~\mathcal{IT})$ be an $\mathcal{ITS}$ and $A,~B\subset X.$ Then
\begin{enumerate}[(i)]
	\item  $\mathcal{I.}e\mathcal{DS}(\emptyset)=\emptyset.$
	\item  If $A\subset B$ then $\mathcal{I.}e\mathcal{DS}(A)\subset \mathcal{I.}e\mathcal{DS}(B)$
	\item if $x\in \mathcal{I.}e\mathcal{DS}(A)$ then $x\in \mathcal{I.}e\mathcal{DS}(A\backslash \{x\})$ . 
	\item $\mathcal{I.}e\mathcal{DS}(A\cap B)\subset \mathcal{I.}e\mathcal{DS}(A)\cap \mathcal{I.}e\mathcal{DS}(B)$.
	\item $\mathcal{I.}e\mathcal{DS}(A\cup B)=\mathcal{I.}e\mathcal{DS}(A)\cup \mathcal{I.}e\mathcal{DS}(B)$ .
\end{enumerate}

\textbf{Proof:}

(i) Suppose that $ \mathcal{I.}e\mathcal{DS}(\emptyset)\neq\emptyset\rightarrow\exists ~x\in \mathcal{I.}e\mathcal{DS}(A)\ni\emptyset\cap(O\backslash \{x\})\neq\emptyset$

$ ~x\in\emptyset$  and $x\not\in\emptyset.$That is contradiction.

$ \mathcal{I.}e\mathcal{DS}(\emptyset)=\emptyset.$

(ii) : Suppose that $A\subset B$. Let $~x\in \mathcal{I.}e\mathcal{DS}(A)\rightarrow\forall O\ni x,A\cap(O\backslash \{x\})\neq\emptyset.$

$\forall O\ni x,B\cap(O\backslash \{x\})\neq\emptyset.$

$x\in \mathcal{I.}e\mathcal{DS}(B)$

$\mathcal{I.}e\mathcal{DS}(A)\subset \mathcal{I.}e\mathcal{DS}(B)\ .$

(iii) Assume that $x\in \mathcal{I.}e\mathcal{DS}(A)\rightarrow\forall O\ni x, A\cap(O\backslash\{x\})\neq\emptyset$

$\forall O\ni x, A\cap(O\cap\{x\}^c)\neq\emptyset$

$\forall O\ni x, A\cap(O\cap\{x\}^c\cap\{x\}^c)\neq\emptyset$

$\forall O\ni x, A\cap(\{x\}^c\cap O\cap\{x\}^c)\neq\emptyset$

$\forall O\ni x, (A\cap(\{x\}^c)\cap (O\cap\{x\}^c)\neq\emptyset$

$\forall O\ni x,~(A\backslash\{x\})\cap(O\backslash\{x\})\neq\emptyset$

$ x\in \mathcal{I.}e\mathcal{DS}(A\backslash\{x\})$

(iv) Since $A\cap B\subset A\cap A\cap B\subset B$

$ \mathcal{I.}e\mathcal{DS}(A\cap B)\subset \mathcal{I.}e\mathcal{DS}(A)\cap \mathcal{I.}e\mathcal{DS}(A\cap B)\subset \mathcal{I.}e\mathcal{DS}(B))$. 

$ \mathcal{I.}e\mathcal{DS}(A\cap B)\subset \mathcal{I.}e\mathcal{DS}(A)\cap \mathcal{I.}e\mathcal{DS}(B))$.

(v) Since $A\subset A\cup B$ and $B\subset A\cup B$ , then $\mathcal{I.}e\mathcal{DS}(A)\subset \mathcal{I.}e\mathcal{DS}(A\cup B)$ and $\mathcal{I.}e\mathcal{DS}(B)\subset \mathcal{I.}e\mathcal{DS}(A\cup B)$ hence $\mathcal{I.}e\mathcal{DS}(A)\cup \mathcal{I.}e\mathcal{DS}(B)\subset \mathcal{I.}e\mathcal{DS}(A\cup B)$ . Conversely,

Suppose that $x\in \mathcal{I.}e\mathcal{DS}(A\cup B)\rightarrow\forall O\ni x, (A\cup B)\cap(O\backslash \{x\})\neq\emptyset.$

$\forall O\ni x,  A\cap(O\backslash\{x\})  \neq\emptyset\cup B\cap(O\backslash \{x\})\neq\emptyset.$

$ x\in \mathcal{I.}e\mathcal{DS}(A)\cup \mathcal{I.}e\mathcal{DS}(B).$Hence, 

$\mathcal{I.}e\mathcal{DS}(A\cup B)=\mathcal{I.}e\mathcal{DS}(A)\cup \mathcal{I.}e\mathcal{DS}(B)$ .

\thm\em Let $(X,~\mathcal{IT})$ be an $\mathcal{ITS}$ and $A,~B\subset X.$ Then
\begin{enumerate}[(i)]
	\item  $\mathcal{I.}e^*\mathcal{DS}(\emptyset)=\emptyset.$
	\item  If $A\subset B$ then $\mathcal{I.}e^*\mathcal{DS}(A)\subset \mathcal{I.}e^*\mathcal{DS}(B)$
	\item if $x\in \mathcal{I.}e^*\mathcal{DS}(A)$ then $x\in \mathcal{I.}e^*\mathcal{DS}(A\backslash \{x\})$ . 
	\item $\mathcal{I.}e^*\mathcal{DS}(A\cap B)\subset \mathcal{I.}e^*\mathcal{DS}(A)\cap \mathcal{I.}e^*\mathcal{DS}(B)$.
	\item $\mathcal{I.}e^*\mathcal{DS}(A\cup B)=\mathcal{I.}e^*\mathcal{DS}(A)\cup \mathcal{I.}e^*\mathcal{DS}(B)$ .
\end{enumerate}

\textbf{Proof:} It follows Theorem \ref{t3.2}

\thm\em Let $(X,~\mathcal{IT})$ be an $\mathcal{ITS}$ and $A,~B\subset X.$ Then
\begin{enumerate}[(i)]
	\item  $\mathcal{I.}a\mathcal{DS}(\emptyset)=\emptyset.$
	\item  If $A\subset B$ then $\mathcal{I.}a\mathcal{DS}(A)\subset \mathcal{I.}a\mathcal{DS}(B)$
	\item if $x\in \mathcal{I.}a\mathcal{DS}(A)$ then $x\in \mathcal{I.}a\mathcal{DS}(A\backslash \{x\})$ . 
	\item $\mathcal{I.}a\mathcal{DS}(A\cap B)\subset \mathcal{I.}a\mathcal{DS}(A)\cap \mathcal{I.}a\mathcal{DS}(B)$.
	\item $\mathcal{I.}a\mathcal{DS}(A\cup B)=\mathcal{I.}a\mathcal{DS}(A)\cup \mathcal{I.}a\mathcal{DS}(B)$ .
\end{enumerate}

\textbf{Proof:} It follows Theorem \ref{t3.2}

\thm\em\label{t3.5} Let $(X,~\mathcal{IT})$ be an $\mathcal{ITS}$ and $A,B\subset X$.Then
\begin{enumerate}[(i)]
	\item $A$ is $\mathcal{I.}e\mathcal{CS}$ iff  $A=\mathcal{I.}e\mathcal{CL}(A)$ .
	\item $ \mathcal{I.}e\mathcal{CL}(\emptyset)=\emptyset$  and $\mathcal{I.}e\mathcal{CL}(X)=X.$
	\item $\mathcal{I.}e\mathcal{CL}(\mathcal{I.}e\mathcal{CL}(A))=\mathcal{I.}e\mathcal{CL}(A)$ .
	\item If $A\subset B$ then $\mathcal{I.}e\mathcal{CL}(A)\subset \mathcal{I.}e\mathcal{CL}(B)$ .
	\item $\mathcal{I.}e\mathcal{CL}(A\cap B)\subset \mathcal{I.}e\mathcal{CL}(A)\cap \mathcal{I.}e\mathcal{CL}(B)$ .
\end{enumerate}
\textbf{Proof.}

(i) Suppose that $A$ is $\mathcal{I.}e\mathcal{CS}$. Since $A\subset A$ and $A\cap A=A\rightarrow \mathcal{I.}e\mathcal{CL}(A)\subset A$, Also $A\subset \mathcal{I.}e\mathcal{CL}(A)\rightarrow A= \mathcal{I.}e\mathcal{CL}(A)$ . Conversely, Let $A=\mathcal{I.}e\mathcal{CL}(A)$ , obviously, $\mathcal{I.}e\mathcal{CL}(A)$ is the smallest $\mathcal{I.}e\mathcal{CS}$. Hence $A$ is $\mathcal{I.}e\mathcal{CS}$.

(ii) Since $X$  and $\emptyset$ are $\mathcal{I.}e\mathcal{CS}$'s, so by (1.) $\mathcal{I.}e\mathcal{CL}(\emptyset)=\emptyset$  and $\mathcal{I.}e\mathcal{CL}(X)=X.$ 

(iii)Since $\mathcal{I.}e\mathcal{CL}(A)$ is the intersection of all $\mathcal{I.}e\mathcal{CS}$'s are $\mathcal{I.}e\mathcal{CS}$s, then $\mathcal{I.}e\mathcal{CL}(\mathcal{I.}e\mathcal{CL}(A))=\mathcal{I.}e\mathcal{CL}(A)$ .

(iv) Consider $A\subset B$. Since $A\subset \mathcal{I.}e\mathcal{CL}(A)$ and $B\subset \mathcal{I.}e\mathcal{CL}(B)$ , so $\mathcal{I.}e\mathcal{CL}(A)\subset \mathcal{I.}e\mathcal{CL}(B)$. 

(v) Since $(A\cap B\subset A\cap A\cap B\subset B)$ then $\mathcal{I.}e\mathcal{CL}(A\cap B)\subset \mathcal{I.}e\mathcal{CL}(A)$ and $\mathcal{I.}e\mathcal{CL}(A\cap B)\subset \mathcal{I.}e\mathcal{CL}(B)\rightarrow \mathcal{I.}e\mathcal{CL}(A\cap B)\subset \mathcal{I.}e\mathcal{CL}(A)\cap \mathcal{I.}e\mathcal{CL}(B)$ .

\thm\em Let $(X,~\mathcal{IT})$ be an $\mathcal{ITS}$ and $A,B\subset X$.Then
\begin{enumerate}[(i)]
	\item $A$ is $\mathcal{I.}e^*\mathcal{CS}$ iff  $A=\mathcal{I.}e^*\mathcal{CL}(A)$ .
	\item $ \mathcal{I.}e^*\mathcal{CL}(\emptyset)=\emptyset$  and $\mathcal{I.}e^*\mathcal{CL}(X)=X.$
	\item $\mathcal{I.}e^*\mathcal{CL}(\mathcal{I.}e^*\mathcal{CL}(A))=\mathcal{I.}e^*\mathcal{CL}(A)$ .
	\item If $A\subset B$ then $\mathcal{I.}e^*\mathcal{CL}(A)\subset \mathcal{I.}e^*\mathcal{CL}(B)$ .
	\item $\mathcal{I.}e^*\mathcal{CL}(A\cap B)\subset \mathcal{I.}e^*\mathcal{CL}(A)\cap \mathcal{I.}e^*\mathcal{CL}(B)$ .
\end{enumerate}

\textbf{Proof:} It follows Theorem \ref{t3.5}

\thm\em Let $(X,~\mathcal{IT})$ be an $\mathcal{ITS}$ and $A,B\subset X$.Then
\begin{enumerate}[(i)]
	\item $A$ is $\mathcal{I.}a\mathcal{CS}$ iff  $A=\mathcal{I.}a\mathcal{CL}(A)$ .
	\item $ \mathcal{I.}a\mathcal{CL}(\emptyset)=\emptyset$  and $\mathcal{I.}a\mathcal{CL}(X)=X.$
	\item $\mathcal{I.}a\mathcal{CL}(\mathcal{I.}a\mathcal{CL}(A))=\mathcal{I.}a\mathcal{CL}(A)$ .
	\item If $A\subset B$ then $\mathcal{I.}a\mathcal{CL}(A)\subset \mathcal{I.}a\mathcal{CL}(B)$ .
	\item $\mathcal{I.}a\mathcal{CL}(A\cap B)\subset \mathcal{I.}a\mathcal{CL}(A)\cap \mathcal{I.}a\mathcal{CL}(B)$ .
\end{enumerate}

\textbf{Proof:} It follows Theorem \ref{t3.5}

\thm\em\label{t3.8} Let $(X,~\mathcal{IT})$ be an $\mathcal{ITS}$ and $A,B\subset X.$ Then
\begin{enumerate}[(i)]
\item $A$ is $\mathcal{I.}e\mathcal{OS}$ iff  $A=\mathcal{I.}e\mathcal{INT}(A)$ . 
\item $ \mathcal{I.}e\mathcal{INT}(X)=X$ and $\mathcal{I.}e\mathcal{INT}(\emptyset)=\emptyset.$
\item $\mathcal{I.}e\mathcal{INT}(\mathcal{I.}e\mathcal{INT}\subset(A))=\mathcal{I.}e\mathcal{INT}(A)$ .
\item If $A\subset B$ then $\mathcal{I.}e\mathcal{INT}(A)\subset \mathcal{I.}e\mathcal{INT}(B)$ . 
\item $\mathcal{I.}e\mathcal{INT}(A\cap B)=\mathcal{I.}e\mathcal{INT}(A)\cap \mathcal{I.}e\mathcal{INT}(B)$ .
\end{enumerate}
\textbf{Proof.}

(i) Suppose that $A$ is $\mathcal{I.}e\mathcal{OS}$. Since $A\subset A$, then $A$ is $\mathcal{I.}e\mathcal{OS}$ containing itself, so $A\subset \mathcal{I.}e\mathcal{INT}(A)$ and $\mathcal{I.}e\mathcal{INT}(A)\subset A$, that implies $A=\mathcal{I.}e\mathcal{INT}(A)$ . Conversely, Let $A=\mathcal{I.}e\mathcal{INT}(A),$ suppose that $A=\mathcal{I.}e\mathcal{INT}(A)$ . Since $\mathcal{I.}e\mathcal{INT}(A)$ is $\mathcal{I.}e\mathcal{OS}$, then $A$ is $\mathcal{I.}e\mathcal{OS}$.

(ii) Since $ X,~~\emptyset$ are $\mathcal{I.}e\mathcal{OS}$'s, by (1), we have $\mathcal{I.}e\mathcal{INT}(X)=X$ and $\mathcal{I.}e\mathcal{INT}(\emptyset)=\emptyset.$

(iii) Since $\mathcal{I.}e\mathcal{INT}(A)$ is $\mathcal{I.}e\mathcal{OS}$. so by (1) $\mathcal{I.}e\mathcal{INT}(\mathcal{I.}e\mathcal{INT}(A))=\mathcal{I.}e\mathcal{INT}(A)$ .

(iv) Suppose that If $A\subset B.$ Let $O_{i}\in \mathcal{I.}e\mathcal{INT}(A)\rightarrow O_{i}\subset A\rightarrow O_{i}\subset B\rightarrow O_{i}\in \mathcal{I.}e\mathcal{INT}(B)$. Therefore $\mathcal{I.}e\mathcal{INT}(A)\subset \mathcal{I.}e\mathcal{INT}(B)$.

(v) Let $O_{i}\in \mathcal{I.}e\mathcal{INT}(A)\cap \mathcal{I.}e\mathcal{INT}(B)\leftrightarrow O_{i}\in \mathcal{I.}e\mathcal{INT}(A)\cap O_{i}\in \mathcal{I.}e\mathcal{INT}(B)$. 

$\leftrightarrow\cup O_{i},O_{i}\subset A,\forall i\cap\cup O_{i},O_{i}\subset B,\forall i.$

$\leftrightarrow\cup O_{i},O_{i}\subset A\cap B,\forall i.$

$\leftrightarrow O_{i}\in \mathcal{I.}e\mathcal{INT}(A\cap B),\forall i.$

\thm\em Let $(X,~\mathcal{IT})$ be an $\mathcal{ITS}$ and $A,B\subset X.$ Then
\begin{enumerate}[(i)]
	\item $A$ is $\mathcal{I.}e^*\mathcal{OS}$ iff  $A=\mathcal{I.}e^*\mathcal{INT}(A)$ . 
	\item $ \mathcal{I.}e^*\mathcal{INT}(X)=X$ and $\mathcal{I.}e^*\mathcal{INT}(\emptyset)=\emptyset.$
	\item $\mathcal{I.}e^*\mathcal{INT}(\mathcal{I.}e^*\mathcal{INT}\subset(A))=\mathcal{I.}e^*\mathcal{INT}(A)$ .
	\item If $A\subset B$ then $\mathcal{I.}e^*\mathcal{INT}(A)\subset \mathcal{I.}e^*\mathcal{INT}(B)$ . 
	\item $\mathcal{I.}e^*\mathcal{INT}(A\cap B)=\mathcal{I.}e^*\mathcal{INT}(A)\cap \mathcal{I.}e^*\mathcal{INT}(B)$ .
\end{enumerate}

\textbf{Proof:} It follows Theorem \ref{t3.8}

\thm\em Let $(X,~\mathcal{IT})$ be an $\mathcal{ITS}$ and $A,B\subset X.$ Then
\begin{enumerate}[(i)]
	\item $A$ is $\mathcal{I.}a\mathcal{OS}$ iff  $A=\mathcal{I.}a\mathcal{INT}(A)$ . 
	\item $ \mathcal{I.}a\mathcal{INT}(X)=X$ and $\mathcal{I.}a\mathcal{INT}(\emptyset)=\emptyset.$
	\item $\mathcal{I.}a\mathcal{INT}(\mathcal{I.}a\mathcal{INT}\subset(A))=\mathcal{I.}a\mathcal{INT}(A)$ .
	\item If $A\subset B$ then $\mathcal{I.}a\mathcal{INT}(A)\subset \mathcal{I.}a\mathcal{INT}(B)$ . 
	\item $\mathcal{I.}a\mathcal{INT}(A\cap B)=\mathcal{I.}a\mathcal{INT}(A)\cap \mathcal{I.}a\mathcal{INT}(B)$ .
\end{enumerate}

\textbf{Proof:} It follows Theorem \ref{t3.8}

\thm\em\label{t3.10} Let $(X,\mathcal{IT})$ be an $\mathcal{ITS}$ and $A,B\subset X.$ Then 
\begin{enumerate}[(i)]
\item $\mathcal{I.}e\mathcal{EP}(X)=\emptyset$ and $\mathcal{I.}e\mathcal{EP}(\emptyset)=X.$
\item $\mathcal{I.}e\mathcal{EP}(A)\subset A^{c}.$
\item $\mathcal{I.}e\mathcal{EP}(A\cup B)=\mathcal{I.}e\mathcal{EP}(A)\cap \mathcal{I.}e\mathcal{EP}(B)$ . 
\item If $A\subset B,$ then $\mathcal{I.}e\mathcal{EP}(B)\subset \mathcal{I.}e\mathcal{EP}(A)$ . 
\item $\mathcal{I.}e\mathcal{EP}(A\cap B)\subset \mathcal{I.}e\mathcal{EP}(A)\cup \mathcal{I.}e\mathcal{EP}(B)$ .
\end{enumerate}
 
\textbf{Proof.}

(i):$\mathcal{I.}e\mathcal{EP}(X)=\mathcal{I.}e\mathcal{INT}(X^{c})=\mathcal{I.}e\mathcal{INT}(\emptyset)=\emptyset$ and $\mathcal{I.}e\mathcal{EP}(\emptyset)=\mathcal{I.}e\mathcal{INT}(\emptyset^{c})=\mathcal{I.}e\mathcal{INT}(X)=X.$

(ii)$\mathcal{I.}e\mathcal{EP}(A)=\mathcal{I.}e\mathcal{INT}(A^{c})\subset A^{c}.$

(iii) $\mathcal{I.}e\mathcal{EP}(A\cup B )=\mathcal{I.}e\mathcal{INT}(A\cup B)^c=\mathcal{I.}e\mathcal{INT}(A^c\cap B^c)=\mathcal{I.}e\mathcal{INT}(A^c)\cap \mathcal{I.}e\mathcal{INT}(B^c)=\mathcal{I.}e\mathcal{EP}(A)\cap \mathcal{I.}e\mathcal{EP}(B)$

(iv) let $A\subset B\rightarrow B^{c}\subset A^{c}\rightarrow \mathcal{I.}e\mathcal{INT}(B^{c})\subset \mathcal{I.}e\mathcal{INT}(A^{c})\rightarrow \mathcal{I.}e\mathcal{EP}(B)\subset \mathcal{I.}e\mathcal{EP}(A)$.

(v) $\mathcal{I.}e\mathcal{EP}(A\cap B)=\mathcal{I.}e\mathcal{INT}(A\cap B)^{c}=\mathcal{I.}e\mathcal{INT}(A^{c}\cup B^{c})\subset \mathcal{I.}e\mathcal{INT}(A^{c})\cup \mathcal{I.}e\mathcal{INT}(B^{c})=\mathcal{I.}e\mathcal{EP}(A)\cup \mathcal{I.}e\mathcal{EP}(B)$.

\thm\em Let $(X,\mathcal{IT})$ be an $\mathcal{ITS}$ and $A,B\subset X.$ Then 
\begin{enumerate}[(i)]
	\item $\mathcal{I.}e^*\mathcal{EP}(X)=\emptyset$ and $\mathcal{I.}e^*\mathcal{EP}(\emptyset)=X.$
	\item $\mathcal{I.}e^*\mathcal{EP}(A)\subset A^{c}.$
	\item $\mathcal{I.}e^*\mathcal{EP}(A\cup B)=\mathcal{I.}e^*\mathcal{EP}(A)\cap \mathcal{I.}e^*\mathcal{EP}(B)$ . 
	\item If $A\subset B,$ then $\mathcal{I.}e^*\mathcal{EP}(B)\subset \mathcal{I.}e^*\mathcal{EP}(A)$ . 
	\item $\mathcal{I.}e^*\mathcal{EP}(A\cap B)\subset \mathcal{I.}e^*\mathcal{EP}(A)\cup \mathcal{I.}e^*\mathcal{EP}(B)$ .
\end{enumerate}
\textbf{Proof:} It follows Theorem \ref{t3.10}

\thm\em Let $(X,\mathcal{IT})$ be an $\mathcal{ITS}$ and $A,B\subset X.$ Then 
\begin{enumerate}[(i)]
	\item $\mathcal{I.}a\mathcal{EP}(X)=\emptyset$ and $\mathcal{I.}a\mathcal{EP}(\emptyset)=X.$
	\item $\mathcal{I.}a\mathcal{EP}(A)\subset A^{c}.$
	\item $\mathcal{I.}a\mathcal{EP}(A\cup B)=\mathcal{I.}a\mathcal{EP}(A)\cap \mathcal{I.}a\mathcal{EP}(B)$ . 
	\item If $A\subset B,$ then $\mathcal{I.}a\mathcal{EP}(B)\subset \mathcal{I.}a\mathcal{EP}(A)$ . 
	\item $\mathcal{I.}a\mathcal{EP}(A\cap B)\subset \mathcal{I.}a\mathcal{EP}(A)\cup \mathcal{I.}a\mathcal{EP}(B)$ .
\end{enumerate}

\textbf{Proof:} It follows Theorem \ref{t3.10}

\thm\em\label{t3.14} Let $(X,~\mathcal{IT})$ be an (ITS) and $A\subset X.$ 
\begin{enumerate}[(i)]
	\item $\mathcal{I.}e\mathcal{B}(X)=\mathcal{I.}e\mathcal{B}(\emptyset)=\emptyset.$
	\item $\mathcal{I.}e\mathcal{B}(A\cap B)=\mathcal{I.}e\mathcal{B}(A)\cup \mathcal{I.}e\mathcal{B}(B)$ .
\end{enumerate}

\textbf{Proof.}

(i) $\mathcal{I.}e\mathcal{B}(X)=X\backslash \mathcal{I.}e\mathcal{INT}(X)\cup \mathcal{I.}e\mathcal{EP}(X)=X\backslash X\cup\emptyset=X\backslash X=\emptyset.$

$ \mathcal{I.}e\mathcal{B}(\emptyset)=X\backslash \mathcal{I.}e\mathcal{INT}(\emptyset)\cup \mathcal{I.}e\mathcal{EP}(\emptyset)=X\backslash \emptyset \cup X=X\backslash X=\emptyset$

(ii) $\mathcal{IB}(A\cap B)=X\backslash \mathcal{I.}e\mathcal{INT}(A\cap B)\cup \mathcal{I.}e\mathcal{EP}(A\cap B)$

$=X\backslash \mathcal{I.}e\mathcal{INT}(A)\cap \mathcal{I.}e\mathcal{INT}(B)\cup \mathcal{I.}e\mathcal{EP}(A\cap B)
$

$
=X\backslash \mathcal{I.}e\mathcal{INT}(A)\cup X\backslash \mathcal{I.}e\mathcal{INT}(B)\cup \mathcal{I.}e\mathcal{EP}(A\cap B)
$

$
=X\backslash \mathcal{I.}e\mathcal{INT}(A)\cup X\backslash \mathcal{I.}e\mathcal{INT}(B)\cup \mathcal{I.}e\mathcal{EP}(A)\cup \mathcal{I.}e\mathcal{EP}(B))
$

$=\mathcal{I.}e\mathcal{B}(A)\cup \mathcal{I.}e\mathcal{B}(B)\ .$

\thm\em Let $(X,~\mathcal{IT})$ be an (ITS) and $A\subset X.$ 
\begin{enumerate}[(i)]
	\item $\mathcal{I.}e^*\mathcal{B}(X)=\mathcal{I.}e^*\mathcal{B}(\emptyset)=\emptyset.$
	\item $\mathcal{I.}e^*\mathcal{B}(A\cap B)=\mathcal{I.}e^*\mathcal{B}(A)\cup \mathcal{I.}e^*\mathcal{B}(B)$ .
\end{enumerate}

\textbf{Proof:} It follows Theorem \ref{t3.14}

\thm\em Let $(X,~\mathcal{IT})$ be an (ITS) and $A\subset X.$ 
\begin{enumerate}[(i)]
	\item $\mathcal{I.}a\mathcal{B}(X)=\mathcal{I.}a\mathcal{B}(\emptyset)=\emptyset.$
	\item $\mathcal{I.}a\mathcal{B}(A\cap B)=\mathcal{I.}a\mathcal{B}(A)\cup \mathcal{I.}a\mathcal{B}(B)$ .
\end{enumerate}

\textbf{Proof:} It follows Theorem \ref{t3.14}

\thm\em\label{t3.17} Let $(X,~\mathcal{IT})$ be an $\mathcal{ITS}$ and $A\subset X.$ Then: 
\begin{enumerate}[(i)]
	\item $A\subset \mathcal{I.}e\mathcal{CL}(A)\rightarrow \mathcal{I.}e\mathcal{DS}(A)\subset \mathcal{I.}e\mathcal{DS}(\mathcal{I.}e\mathcal{CL}(A))$.
	\item $\mathcal{I.}e\mathcal{INT}(A)\subset A\rightarrow \mathcal{I.}e\mathcal{DS}(\mathcal{I.}e\mathcal{INT}(A))\subset \mathcal{I.}e\mathcal{DS}(A)$ .
	\item If $A$ is $\mathcal{ICS}$, then $\mathcal{I.}e\mathcal{DS}(A)\subset A.$
	\item $\mathcal{I.}e\mathcal{CL}(A)=A\cup \mathcal{I.}e\mathcal{DS}(A)$ .
	\item  $\mathcal{IB}(A)=\mathcal{I.}e\mathcal{CL}(A)\backslash \mathcal{I.}e\mathcal{INT}(A)$ .
	\item  $\mathcal{I.}e\mathcal{CL}(A)=\mathcal{IB}(A)\cup \mathcal{I.}e\mathcal{INT}(A)$ .
	\item  $\mathcal{IB}(A)\subset \mathcal{I.}e\mathcal{CL}(A)$ .
	\item  $\mathcal{I.}e\mathcal{INT}(A)\cap \mathcal{IB}(A)=\emptyset.$
\end{enumerate}

\textbf{Proof.}

(i) Let $A\subset \mathcal{I.}e\mathcal{CL}(A)$. By (ii) $\mathcal{I.}e\mathcal{DS}(A)\subset \mathcal{I.}e\mathcal{DS}(\mathcal{I.}e\mathcal{CL}(A))$.

(ii) Let $\mathcal{I.}e\mathcal{INT}(A)\subset A.$ By (ii) $\mathcal{I.}e\mathcal{DS}(\mathcal{I.}e\mathcal{INT}(A))\subset \mathcal{I.}e\mathcal{DS}(A)$ .

(iii) Let $A$ be a $\mathcal{ICS}$ and $x\in \mathcal{I.}e\mathcal{DS}(A)$ , then $\forall O\ni x, A\cap(O-(x))\neq\emptyset$ Hence $x\in A$ and $\mathcal{I.}e\mathcal{DS}(A)\subset A.$

(iv) Since $A\subset \mathcal{I.}e\mathcal{CL}(A)$ and $\mathcal{I.}e\mathcal{DS}(A)\subset \mathcal{I.}e\mathcal{DS}(\mathcal{I.}e\mathcal{CL}(A))\subset \mathcal{I.}e\mathcal{CL}(A)$ . we have $A\cup \mathcal{I.}e\mathcal{DS}(A)\subset \mathcal{I.}e\mathcal{CL}(A)$ . Another direction, To show that $\mathcal{I.}e\mathcal{CL}(A)\subset A\cup \mathcal{I.}e\mathcal{DS}(A)$ .

Let $x\in \mathcal{I.}e\mathcal{CL}(A)$ , but $A\subset \mathcal{I.}e\mathcal{CL}(A)$ , then $x\in A$ or $x\not\in A.$

(Probability 1) If $x\in A$, then $x\in A\cup \mathcal{I.}e\mathcal{DS}(A)$ .

(Probability 2) If $x\not\in A,$ Let $x\not\in \mathcal{I.}e\mathcal{DS}(A)\rightarrow\exists O\ni x,A\cap(O\backslash \{x\})=\emptyset$, but $x\not\in A$, that is contradiction, therefore $x\in \mathcal{I.}e\mathcal{DS}(A)$ and $x\in A\cup \mathcal{I.}e\mathcal{DS}(A)$ .

So $\mathcal{I.}e\mathcal{CL}(A)=A\cup \mathcal{I.}e\mathcal{DS}(A)$ .

(v) By definition: $\mathcal{IB}(A)=X\backslash \mathcal{I.}e\mathcal{INT}(A)\cup iep(A)$

$
=X\backslash \mathcal{I.}e\mathcal{INT}(A)\cap X\backslash iep(A)
$

$
=X\backslash \mathcal{I.}e\mathcal{INT}(A)\cap \mathcal{I.}e\mathcal{CL}\backslash (A)
$

Since $\mathcal{I.}e\mathcal{INT}(A)\subset \mathcal{I.}e\mathcal{CL}(A)\subset X\rightarrow \mathcal{I.}e\mathcal{CL}(A)\cap(X\backslash \mathcal{I.}e\mathcal{INT}(A)=\mathcal{I.}e\mathcal{CL}(A)\backslash \mathcal{I.}e\mathcal{INT}(A)$ . Then we have
$
\mathcal{IB}(A)=\mathcal{I.}e\mathcal{CL}(A)\backslash \mathcal{I.}e\mathcal{INT}(A)\ .
$

(vi) By (i) $\mathcal{IB}(A)=\mathcal{I.}e\mathcal{CL}(A)\backslash \mathcal{I.}e\mathcal{INT}(A)$

$
\rightarrow \mathcal{IB}(A)\cup \mathcal{I.}e\mathcal{INT}(A)=ic(A)\backslash \mathcal{I.}e\mathcal{INT}(A)\cup \mathcal{I.}e\mathcal{INT}(A)=\mathcal{I.}e\mathcal{CL}(A)\ .$

(vii) By (ii) it is clear that $\mathcal{IB}(A)\subset \mathcal{I.}e\mathcal{CL}(A)$ .

(viii) $\mathcal{I.}e\mathcal{INT}(A)\cap \mathcal{IB}(A)=\mathcal{I.}e\mathcal{INT}(A)\cap \mathcal{I.}e\mathcal{CL}(A)\backslash \mathcal{I.}e\mathcal{INT}(A)=\emptyset.$

\thm\em Let $(X,~\mathcal{IT})$ be an $\mathcal{ITS}$ and $A\subset X.$ Then: 
\begin{enumerate}[(i)]
	\item $A\subset \mathcal{I.}e^*\mathcal{CL}(A)\rightarrow \mathcal{I.}e^*\mathcal{DS}(A)\subset \mathcal{I.}e^*\mathcal{DS}(\mathcal{I.}e^*\mathcal{CL}(A))$.
	\item $\mathcal{I.}e^*\mathcal{INT}(A)\subset A\rightarrow \mathcal{I.}e^*\mathcal{DS}(\mathcal{I.}e^*\mathcal{INT}(A))\subset \mathcal{I.}e^*\mathcal{DS}(A)$ .
	\item If $A$ is $\mathcal{ICS}$, then $\mathcal{I.}e^*\mathcal{DS}(A)\subset A.$
	\item $\mathcal{I.}e^*\mathcal{CL}(A)=A\cup \mathcal{I.}e^*\mathcal{DS}(A)$ .
	\item  $\mathcal{IB}(A)=\mathcal{I.}e^*\mathcal{CL}(A)\backslash \mathcal{I.}e^*\mathcal{INT}(A)$ .
	\item  $\mathcal{I.}e^*\mathcal{CL}(A)=\mathcal{IB}(A)\cup \mathcal{I.}e^*\mathcal{INT}(A)$ .
	\item  $\mathcal{IB}(A)\subset \mathcal{I.}e^*\mathcal{CL}(A)$ .
	\item  $\mathcal{I.}e^*\mathcal{INT}(A)\cap \mathcal{IB}(A)=\emptyset.$
\end{enumerate}

\textbf{Proof:} It follows Theorem \ref{t3.17}

\thm\em Let $(X,~\mathcal{IT})$ be an $\mathcal{ITS}$ and $A\subset X.$ Then: 
\begin{enumerate}[(i)]
	\item $A\subset \mathcal{I.}a\mathcal{CL}(A)\rightarrow \mathcal{I.}a\mathcal{DS}(A)\subset \mathcal{I.}a\mathcal{DS}(\mathcal{I.}a\mathcal{CL}(A))$.
	\item $\mathcal{I.}a\mathcal{INT}(A)\subset A\rightarrow \mathcal{I.}a\mathcal{DS}(\mathcal{I.}a\mathcal{INT}(A))\subset \mathcal{I.}a\mathcal{DS}(A)$ .
	\item If $A$ is $\mathcal{ICS}$, then $\mathcal{I.}a\mathcal{DS}(A)\subset A.$
	\item $\mathcal{I.}a\mathcal{CL}(A)=A\cup \mathcal{I.}a\mathcal{DS}(A)$ .
	\item  $\mathcal{IB}(A)=\mathcal{I.}a\mathcal{CL}(A)\backslash \mathcal{I.}a\mathcal{INT}(A)$ .
	\item  $\mathcal{I.}a\mathcal{CL}(A)=\mathcal{IB}(A)\cup \mathcal{I.}a\mathcal{INT}(A)$ .
	\item  $\mathcal{IB}(A)\subset \mathcal{I.}a\mathcal{CL}(A)$ .
	\item  $\mathcal{I.}a\mathcal{INT}(A)\cap \mathcal{IB}(A)=\emptyset.$
\end{enumerate}

\textbf{Proof:} It follows Theorem \ref{t3.17}

\begin{thm} \label{l1}
	\emph{  Let $(X,~\mathcal{IT})$ be an $\mathcal{ITS}$ and $A\subset X.$ Then: 
		\begin{enumerate}
			\item[(i)] $\mathcal{I.}\delta\mathcal{PCL}(A)\supseteq A\cup \mathcal{I.CL}(\mathcal{I.INT_\delta}(A))$ and $\mathcal{I.}\delta\mathcal{PNT}(A)\subseteq A\cap \mathcal{I.INT}(\mathcal{I.CL_\delta}(A))$
			\item[(ii)] $\mathcal{I.}\delta\mathcal{SCL}(A) \supseteq A\cup \mathcal{I.INT}(\mathcal{I.CL_\delta}(A))$ and $\mathcal{I.}\delta\mathcal{SINT}(A) \subseteq A\cap \mathcal{I.CL}(\mathcal{I.INT_\delta}(A))$.
	\end{enumerate}}
\end{thm}
{\NI\bf Proof.} We will prove only the first statement of (i) and the others is similar. Since $\mathcal{I.}\delta\mathcal{PCL}(A)$ is $\mathcal{I.}\delta\mathcal{PCS}$, we have $\mathcal{I.CL}(\mathcal{I.INT_\delta}(A))\subseteq \mathcal{I.CL}(\mathcal{I.}\delta\mathcal{PCL}(A))\subseteq \mathcal{I.}\delta\mathcal{PCL}(A)$. Thus $A\cup \mathcal{I.CL}(\mathcal{I.INT_\delta}
(A))\subseteq \mathcal{I.}\delta\mathcal{PCL}(A)$. \hfill\eb

\begin{prop}
	\emph{ Let $(X,~\mathcal{IT})$ be an $\mathcal{ITS}$ and $A\subset X.$ Then:
		\begin{enumerate}
			\item[(i)] If $A$ is an $\mathcal{I.}e\mathcal{OS}$ and $\mathcal{I.INT_\delta}(A)=\phi$, then $A$ is an $\mathcal{I.}\delta\mathcal{POS}$.
			\item[(ii)] If $A$ is an $\mathcal{I.}e\mathcal{OS}$ and $\mathcal{I.CL_\delta}(A)=\phi$, then $A$ is an $\mathcal{I.}\delta\mathcal{SOS}$.
			\item[(iii)] If $A$ is an $\mathcal{I.}e\mathcal{OS}$ and $\mathcal{I.}\delta\mathcal{CS}$, then $A$ is an $\mathcal{I.}\delta\mathcal{SOS}$.
			\item[(iv)] If $A$ is an $\mathcal{I.}\delta\mathcal{SOS}$ and $\mathcal{I.}\delta\mathcal{CS}$, then $A$ is an $\mathcal{I.}e\mathcal{OS}$.
	\end{enumerate}}
\end{prop}
{\NI\bf Proof.} (i) Let $A$ be an $\mathcal{I.}e\mathcal{OS}$, that is $A\subseteq \mathcal{I.CL}\mathcal{I.INT_\delta}(A)\cup \mathcal{I.INT}\mathcal{I.CL_\delta}(A)=\phi\cup \mathcal{I.INT}\mathcal{I.CL_\delta}(A)=\mathcal{I.INT}\mathcal{I.CL_\delta}(A)$. Hence $A$ is an $\mathcal{IPOS}$.

(ii) Follows from (i).

(iii) Let $A$ be an $\mathcal{I.}e\mathcal{OS}$ and $\mathcal{I.}\delta\mathcal{CS}$, that is $A\subseteq \mathcal{I.CL}\mathcal{I.INT_\delta}(A)\cup \mathcal{I.INT}\mathcal{I.CL_\delta}(A)=\mathcal{I.CL}\mathcal{I.INT_\delta}(A)\cup \mathcal{I.INT}(A)=cl\mathcal{I.INT_\delta}(A)$. Hence $A$ is an $\mathcal{I.}\delta\mathcal{SOS}$.

(iv) Let $A$ be an $\mathcal{I.}\delta\mathcal{SOS}$ and $\mathcal{I.}\delta\mathcal{CS}$, that is $A\subseteq \mathcal{I.CL}\mathcal{I.INT_\delta}(A)\subseteq \mathcal{I.CL}\mathcal{I.INT_\delta}(A)\vee \mathcal{I.INT}\mathcal{I.CL_\delta}(A)$. Hence $A$ is an $\mathcal{I.}e\mathcal{OS}$.
\hfill\eb

\begin{thm}
	\emph{ Let $(X,~\mathcal{IT})$ be an $\mathcal{ITS}$ and $A\subset X.$ Then, $A$ is an $\mathcal{I.}e\mathcal{OS}$ if and only if $A=\mathcal{I.}\delta\mathcal{PINT}(A)\cup \mathcal{I.}\delta\mathcal{SINT}(A)$.}
\end{thm}
{\NI\bf Proof.} Let $A$ be an $\mathcal{I.}e\mathcal{OS}$. Then $A\subseteq \mathcal{I.CL}(\mathcal{I.INT_\delta}(A))\cup \mathcal{I.INT}(\mathcal{I.CL_\delta}(A))$. By Theorem \ref{l1}, we have $\mathcal{I.}\delta\mathcal{PINT}(A)\cup \mathcal{I.}\delta\mathcal{SINT}(A)=(A\cap \mathcal{I.INT}(\mathcal{I.CL_\delta}(A)))\cup  (A\cap \mathcal{I.CL}(\mathcal{I.INT_\delta}(A)))= A\cap (\mathcal{I.INT}(\mathcal{I.CL_\delta}(A))\cup  \mathcal{I.CL}(\mathcal{I.INT_\delta}(A)))=A $.

Conversely, if $A=\mathcal{I.}\delta\mathcal{PINT}(A)\cup \mathcal{I.}\delta\mathcal{SINT}(A)$ then, by Theorem \ref{l1} $A=\mathcal{I.}\delta\mathcal{PINT}(A)\cup \mathcal{I.}\delta\mathcal{SINT}(A)=(A\cap \mathcal{I.INT}(\mathcal{I.CL_\delta}(A)))\cup  (A\cap \mathcal{I.CL}(\mathcal{I.INT_\delta}(A)))= A\cap (\mathcal{I.INT}(\mathcal{I.CL_\delta}(A))\cup  \mathcal{I.CL}(\mathcal{I.INT_\delta}(A)))\subseteq \mathcal{I.INT}(\mathcal{I.CL_\delta}(A))\cup  \mathcal{I.CL}(\mathcal{I.INT_\delta}(A))$ and hence $A$ is an $\mathcal{I.}e\mathcal{OS}$. \hfill\eb

\begin{prop}
	\emph{ Let $(X,~\mathcal{IT})$ be an $\mathcal{ITS}$ and $A\subset X.$ Then:
		\begin{enumerate}
			\item[(i)] $\mathcal{I.}e\mathcal{CL}(\overline{A})=\overline{\mathcal{I.}e\mathcal{INT}(A)}$, $\mathcal{I.}e\mathcal{INT}(\overline{A})=\overline{\mathcal{I.}e\mathcal{CL}(A)}$.
			\item[(ii)] $\mathcal{I.}e\mathcal{CL}(A \vee B)\geq \mathcal{I.}e\mathcal{CL}(A)\vee \mathcal{I.}e\mathcal{CL}(B)$, $\mathcal{I.}e\mathcal{INT}(A \vee B)\geq \mathcal{I.}e\mathcal{INT}(A)\vee \mathcal{I.}e\mathcal{INT}(B)$.
			\item[(iii)] $\mathcal{I.}e\mathcal{CL}(A \wedge B)\subseteq \mathcal{I.}e\mathcal{CL}(A)\wedge \mathcal{I.}e\mathcal{CL}(B)$, $\mathcal{I.}e\mathcal{INT}(A \wedge B)\subseteq \mathcal{I.}e\mathcal{INT}(A)\wedge \mathcal{I.}e\mathcal{INT}(B)$.
	\end{enumerate}}
\end{prop}

\begin{prop} \label{p1}
	\emph{ Let $(X,~\mathcal{IT})$ be an $\mathcal{ITS}$ and $A\subset X.$ Then:
		\begin{enumerate}
			\item[(i)] $\mathcal{I.}e\mathcal{CL}(A)\geq \mathcal{I.CL}\mathcal{I.INT_\delta}(A) \wedge \mathcal{I.INT}\mathcal{I.CL_\delta}(A)$.
			\item[(ii)] $\mathcal{I.}e\mathcal{INT}(A)\subseteq \mathcal{I.CL}\mathcal{I.INT_\delta}(A)\vee \mathcal{I.INT}\mathcal{I.CL_\delta}(A)$.
	\end{enumerate}}
\end{prop}
{\NI\bf Proof.} (i) $\mathcal{I.}e\mathcal{CL}(A)$ is an $\mathcal{I.}e\mathcal{CS}$ and $A\subseteq \mathcal{I.}e\mathcal{CL}(A)$, then $\mathcal{I.}e\mathcal{CL}(A)\geq \mathcal{I.CL}\mathcal{I.INT_\delta}\mathcal{I.}e\mathcal{CL}(A)\wedge \mathcal{I.INT}\mathcal{I.CL_\delta}\mathcal{I.}e\mathcal{CL}(A)\geq \mathcal{I.CL}\mathcal{I.INT_\delta}(A)\wedge \mathcal{I.INT}\mathcal{I.CL_\delta}(A)$.

(ii) Follows from (i) by taking the complementation.
\hfill\eb

\begin{thm}
	\emph{ Let $(X,~\mathcal{IT})$ be an $\mathcal{ITS}$ and $A\subset X.$ Then: $\mathcal{I.}e\mathcal{CL}(A)=\mathcal{I.}\delta\mathcal{PCL}(A)\cap \mathcal{I.}\delta\mathcal{SCL}(A)$.}
\end{thm}
{\NI\bf Proof.} It is obvious that, $\mathcal{I.}e\mathcal{CL}(A) \subseteq \mathcal{I.}\delta\mathcal{PCL}(A)\cap \mathcal{I.}\delta\mathcal{SCL}(A)$. Conversely, from Definition we have $\mathcal{I.}e\mathcal{CL}(A)\supseteq \mathcal{I.CL}(\mathcal{I.INT_\delta}(\mathcal{I.}e\mathcal{CL}(A)))\cap \mathcal{I.INT}(\mathcal{I.CL_\delta}(\mathcal{I.}e\mathcal{CL}(A)))\supseteq \mathcal{I.CL}(\mathcal{I.INT_\delta}(A))\cap \mathcal{I.INT}(\mathcal{I.CL_\delta}(A))$. Since $\mathcal{I.}e\mathcal{CL}(A)$ is $\mathcal{I.}e\mathcal{OS}$, by Theorem \ref{l1}, we have $\mathcal{I.}\delta\mathcal{PCL}(A)\cap \mathcal{I.}\delta\mathcal{SCL}(A)=(A\cup \mathcal{I.CL}
(\mathcal{I.INT_\delta}(A))) \cup (A\cup \mathcal{I.INT}(\mathcal{I.CL_\delta}(A)))= A\cup (\mathcal{I.CL}(\mathcal{I.INT_\delta}(A)) \cap \mathcal{I.INT}(\mathcal{I.CL_\delta}(A)))=A\subseteq \mathcal{I.}e\mathcal{CL}(A)$. \hfill\eb

\lem\em The following hold for a subset $A$  of a space $X$: 
\begin{enumerate}[(1)]
\item  $\mathcal{I.}\delta\mathcal{SINT}(A)=A\cap \mathcal{I.CL}(\mathcal{I.INT_\delta}(A))$ and 
$\mathcal{I.}\delta\mathcal{SCL}(A)=A\cup \mathcal{I.INT}(\mathcal{I.CL_\delta}(A))$
\item  $\mathcal{I.}\delta\mathcal{PCL}(A)=A\cup \mathcal{I.CL}(\mathcal{I.INT_\delta}(A))$
\item $\mathcal{I.}\delta\mathcal{SCL}(\mathcal{I.}\delta\mathcal{SINT}(A))=\mathcal{I.}\delta\mathcal{SINT}(A)\cup \mathcal{I.INT}(\mathcal{I.CL}(\mathcal{I.INT_\delta}(A)))$  and $ \mathcal{I.}\delta\mathcal{SINT}(\mathcal{I.}\delta\mathcal{SCL}(A))=\mathcal{I.}\delta\mathcal{SCL}(A)\cap \mathcal{I.CL}(\mathcal{I.INT}(\mathcal{I.CL_\delta}(A)))\ $
\item $\mathcal{I.CL_\delta}(\mathcal{I.}\delta\mathcal{SINT}(A))=\mathcal{I.CL}(\mathcal{I.INT_\delta}(A))$
\item $\mathcal{I.}\delta\mathcal{SCL}(\mathcal{I.INT_\delta}(A))=\mathcal{I.INT}(\mathcal{I.CL}(\mathcal{I.INT_\delta}(A)))$
\end{enumerate}







%
%

\lem\em  The following hold for a subset $A$  of a space $X$: 
\begin{enumerate}[(1)]
	\item $\mathcal{I.}e^*\mathcal{CL}( A)$  is $\mathcal{I.}e^*\mathcal{OS}$
	\item $X\backslash \mathcal{I.}e^*\mathcal{CL}(A)=\mathcal{I.}e^*\mathcal{INT}(X\backslash A)$
\end{enumerate}

\thm\em  The following hold for a subset $A$  of a space $X$: 
\begin{enumerate}[(i)]
	\item $A$  is $\mathcal{I.}e^*\mathcal{OS}$ if and only if $A=A\cap \mathcal{I.CL}(\mathcal{I.INT}(\mathcal{I.CL_\delta}(A)))$ 
	\item $A$  is $\mathcal{I.}e^*\mathcal{CS}$ if and only if $A=A\cup \mathcal{I.INT}(\mathcal{I.CL}(\mathcal{I.INT_\delta}(A)))$
	\item $\mathcal{I.}e^*\mathcal{CL}(A)=A\cup \mathcal{I.INT}(\mathcal{I.CL}(\mathcal{I.INT_\delta}(A)))$
	\item $\mathcal{I.}e^*\mathcal{INT}(A)=A\cap \mathcal{I.CL}(\mathcal{I.INT}(\mathcal{I.CL_\delta}(A)))$
\end{enumerate}

\textbf{Proof.} (i) : Let $A$ be $\mathcal{I.}a\mathcal{OS}$. Then $A\subset \mathcal{I.CL}(\mathcal{I.INT}(\mathcal{I.CL_\delta}(A)))$ . We obtain $A\subset A\cap \mathcal{I.CL}(\mathcal{I.INT}(\mathcal{I.CL_\delta}(A)))$ . Conversely, let $A=A\cap \mathcal{I.CL}(\mathcal{I.INT}(\mathcal{I.CL_\delta}(A)))$ . We have $A=A\cap \mathcal{I.CL}(\mathcal{I.INT}(\mathcal{I.CL_\delta}(A)))\subset \mathcal{I.CL}(\mathcal{I.INT}(\mathcal{I.CL_\delta}(A)))$ and hence, $A$ is $\mathcal{I.}e^*\mathcal{OS}$.

(iii):Since $\mathcal{I.}e^*\mathcal{CL}(A)$ is $\mathcal{I.}e^*\mathcal{CS}$,  $\mathcal{I.INT}(\mathcal{I.CL}(\mathcal{I.INT_\delta}(A)))\subset \mathcal{I.INT}(\mathcal{I.CL}(\mathcal{I.INT_\delta}(\mathcal{I.}e^*\mathcal{CL}(A)))) \subset \mathcal{I.}e^*\mathcal{CL}(A)$ . Hence, $A\cup \mathcal{I.INT}(\mathcal{I.CL}(\mathcal{I.INT_\delta}(A)))\subset \mathcal{I.}e^*\mathcal{CL}(A)$ .

Conversely, since:
 $\mathcal{I.INT}(\mathcal{I.CL}(\mathcal{I.INT_\delta}(A\cup \mathcal{I.INT}(\mathcal{I.CL}(\mathcal{I.INT_\delta}(A))))))=$
$$
=\mathcal{I.INT}(\mathcal{I.CL}(\mathcal{I.INT_\delta}(A\cup\mathcal{I.INT_\delta}(\mathcal{I.CL_\delta}(\mathcal{I.INT_\delta}(A))))))
$$
$$
=\mathcal{I.INT}(\mathcal{I.CL}(\mathcal{I.INT_\delta}(A)\cup\mathcal{I.INT_\delta}(\mathcal{I.INT_\delta}(\mathcal{I.CL_\delta}(\mathcal{I.INT_\delta}(A))))))
$$
$$
=\mathcal{I.INT}(\mathcal{I.CL}(\mathcal{I.INT_\delta}(A)\cup\mathcal{I.INT_\delta}(\mathcal{I.CL_\delta}(\mathcal{I.INT_\delta}(A)))))
$$
$$
=\mathcal{I.INT}(\mathcal{I.CL}(\mathcal{I.INT_\delta}(\mathcal{I.CL_\delta}(\mathcal{I.INT_\delta}(A)))))
$$
$$
=\mathcal{I.INT}(\mathcal{I.CL}(\mathcal{I.INT_\delta}(A)))\subset A\cup \mathcal{I.INT}(\mathcal{I.CL}(\mathcal{I.INT_\delta}(A)))\ ,
$$
then $A\cup \mathcal{I.INT}(\mathcal{I.CL}(\mathcal{I.INT_\delta}(A)))$ is $\mathcal{I.}e^*\mathcal{CS}$ containing $A$ and hence:
$$
\mathcal{I.}e^*\mathcal{CL}(A)\subset A\cup \mathcal{I.INT}(\mathcal{I.CL}(\mathcal{I.INT_\delta}(A)))\ .
$$
Thus, we obtain $\mathcal{I.}e^*\mathcal{CL}(A)=A\cup \mathcal{I.INT}(\mathcal{I.CL}(\mathcal{I.INT_\delta}(A)))$ .

(ii) follows from (i) and (iv) follows from (iii).

%
%
%
%
%
%
%
%


\thm\em Let $N$  be a subset of an $\mathcal{ITS}$ $X$. The following are equivalent:

\begin{enumerate}[(i)]
	\item $N$  is $\mathcal{I.}\mathcal{ROS}$,
	
	\item $N$  is $\mathcal{I.}a\mathcal{OS}$ and $\mathcal{I.}e^*\mathcal{CS}$,
	
	\item$N$  is $\mathcal{I.}\delta\mathcal{POS}$ and $\mathcal{I.}\delta\mathcal{SCS}$.
\end{enumerate}

\textbf{Proof.} (i) $\Rightarrow$ (ii) : Obvious.

(ii) $\Rightarrow$(i) : Let $N$ be $\mathcal{I.}a\mathcal{OS}$ and $\mathcal{I.}e^*\mathcal{CS}$. We have $N\subset \mathcal{I.INT}(\mathcal{I.CL}(\mathcal{I.INT_\delta}(N)))$ and  $\mathcal{I.INT}(\mathcal{I.CL}(\mathcal{I.INT_\delta}(N)))\subset N$ and hence $N=\mathcal{I.INT}(\mathcal{I.CL}(\mathcal{I.INT_\delta}(N)))$ . Thus, $N$ is $\mathcal{I.}\mathcal{ROS}$.

(i) $\Leftrightarrow$ (iii) : Let $N$ be $\mathcal{I.}\delta\mathcal{POS}$ and $\mathcal{I.}\delta\mathcal{SCS}$. Then $N\subset \mathcal{I.INT}(\mathcal{I.CL_\delta}(N))$ and  $\mathcal{I.INT}(\mathcal{I.CL_\delta}(N))\subset N$. Thus, $N=\mathcal{I.INT}(\mathcal{I.CL_\delta}(N))=\mathcal{I.INT}(\mathcal{I.CL}(N))$ and hence $N$ is $\mathcal{I.}\mathcal{ROS}$. The converse is similar. 

\thm\em  Let $N$  be a subset of an $\mathcal{ITS}$ $X$. The following are equivalent:
\begin{enumerate}[(i)]
	\item $N$  is $\mathcal{I.}\delta\mathcal{SOS}$,
	
	\item $N$  is $\mathcal{I.}e^*\mathcal{OS}$ and $\mathcal{I.INT_\delta}(\mathcal{I.}\delta\mathcal{FR}(N))=\emptyset.$
\end{enumerate}
\textbf{Proof.} (i) $\Rightarrow$ (ii) : Let $N$ be $\mathcal{I.}\delta\mathcal{SOS}$. We have  $\mathcal{I.INT}(\mathcal{I.CL_\delta}(N))\subset \mathcal{I.CL_\delta}(N)\subset \mathcal{I.CL}(\mathcal{I.INT_\delta}(N))$ . Since 
 $\mathcal{I.INT_\delta}(\mathcal{I.}\delta\mathcal{FR}(N))=\mathcal{I.INT_\delta}(\mathcal{I.CL_\delta}(N)\cap(X\backslash \mathcal{I.INT_\delta}(N)))=\mathcal{I.INT_\delta}(\mathcal{I.CL_\delta}(N))\backslash \mathcal{I.CL}(\mathcal{I.INT_\delta}(N))$ ,
then  $\mathcal{I.INT_\delta}(\mathcal{I.}\delta\mathcal{FR}(N))=\emptyset.$

(ii) $\Rightarrow$ (i) : Let $N$ be $\mathcal{I.}e^*\mathcal{OS}$ and  $\mathcal{I.INT_\delta}(\mathcal{I.}\delta\mathcal{FR}(N))=\emptyset$. Then $N\subset \mathcal{I.CL}(\mathcal{I.INT}(\mathcal{I.CL_\delta}(N)))\subset \mathcal{I.CL}(\mathcal{I.INT_\delta}(N))$ . Thus, $N$ is $\mathcal{I.}\delta\mathcal{SOS}$. 



\thm\em  Let $X$  be a topological space. Then $I.aO(X)=\mathcal{I.}\delta\mathcal{SOS}(X)\cap\mathcal{I.}\delta\mathcal{POS}(X)$ .

\textbf{Proof.} Let $N\in \mathcal{I.}a\mathcal{O}(X)$ . Then $N\in\mathcal{I.}\delta\mathcal{SO}(X)$ and $N\in\mathcal{I.}\delta\mathcal{PO}(X)$ . Thus, $\mathcal{I.}a\mathcal{O}(X)\subset\mathcal{I.}\delta\mathcal{SO}(X)\cap\mathcal{I.}\delta\mathcal{PO}(X)$ .

Conversely, let $N\in\mathcal{I.}$ $\delta\mathcal{SO}(X)\cap\mathcal{I.}\delta\mathcal{PO}(X)$ . Then $N\in\mathcal{I.}\delta\mathcal{SO}(X)$ and $N\in\mathcal{I.}\delta\mathcal{PO}(X)$ . Since $N\in\mathcal{I.}\delta\mathcal{SO}(X)$ , then ,  $\mathcal{I.INT_\delta}(\mathcal{I.}\delta\mathcal{FR}(N))=\emptyset$. Since $\mathcal{I.INT_\delta}(\mathcal{I.}\delta\mathcal{FR}(N))=\mathcal{I.INT_\delta}(\mathcal{I.CL_\delta}(N)\cap(X\backslash \mathcal{I.INT_\delta}(N)))=\mathcal{I.INT_\delta}(\mathcal{I.CL_\delta}(N))\backslash \mathcal{I.CL_\delta}(\mathcal{I.INT_\delta}(N))$, 
then  $\mathcal{I.INT}(\mathcal{I.CL_\delta}(N))\subset \mathcal{I.CL}\mathcal{I.INT_\delta}(N))$ . Since $N\in\mathcal{I.}\delta\mathcal{PO}(X)$ , we have
$N \subset \mathcal{I.INT}(\mathcal{I.CL}_{\delta}(N)) \subset \mathcal{I.INT}(\mathcal{I.CL}(\mathcal{I.INT}_{\delta}(N)))$.
Thus, $N \in \mathcal{I}.aO(X)$.


\section{Interrelations}
\ex\em Let $X$ be a set $X=\{a,b,c,d\},\mathcal{IT}=\{\phi,X,\{a\},\{b\},\{a,c\}\}$ Then 
\begin{enumerate}[(i)]
\item $\{a\}$ is $\mathcal{IOS}$ and $\mathcal{ISOS}$ but not $\mathcal{I.}\delta\mathcal{OS}$ and $\mathcal{I.}\delta\mathcal{SOS}$ 
\item $\{a,b\}$ is $\mathcal{I.}\delta\mathcal{POS}$ but not $\mathcal{IOS}$
\item $\{c\}$ is $\mathcal{I.}\delta\mathcal{POS}$ and $\mathcal{I.}e\mathcal{OS}$ but not $\mathcal{IPOS}$ and $\mathcal{I.}\delta\mathcal{SOS}$
\item $\{b,d\}$ is $\mathcal{I.}e\mathcal{OS}$ but not $\mathcal{I.}\delta\mathcal{POS}$
\item $\{a,b,c\}$ is $\mathcal{I.}e\mathcal{OS}$ but not $\mathcal{IOS}$ 
\item $\{c,d\}$ is $\mathcal{I.}e^*\mathcal{OS}$ but not $\mathcal{I.}e\mathcal{OS}$ 
\item  $\{b,c,d\}$  is $\mathcal{I.}e^*\mathcal{OS}$ but not $\mathcal{I.}\beta\mathcal{OS}$
\end{enumerate}

\ex\em Let $X$ be a set $X=\{a,b,c,d\},\mathcal{IT}=\{\phi,X,\{b\},\{c\}.\{b,c,d\}\}$ Then  $\{b,d\}$ is $\mathcal{I.}\delta\mathcal{SOS}$ but not $\mathcal{IOS}$ 

\includegraphics[scale=0.7]{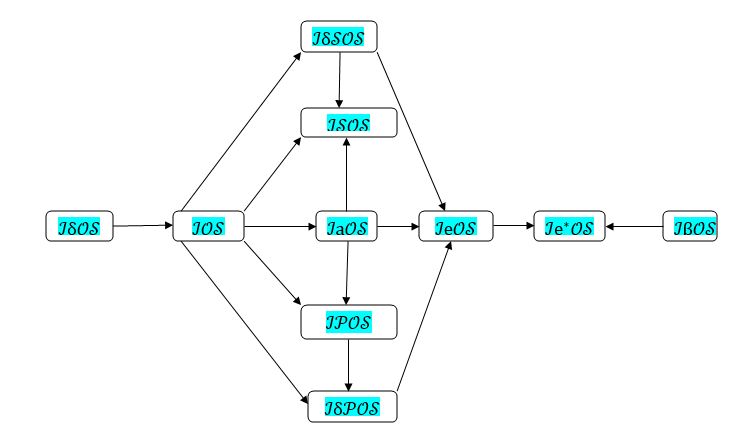}


\begin{thebibliography}{99}
\bibitem{ade} Adel. M. Al. Odhari, \textit{On Infra topological spaces}, Internaational Jounrnal of Mathematical Archive, \textbf{6}(11) 2015, 179-184.

\bibitem{bro} R. Brown, Topology: A geometric Account of General Topology, Homotopy Types and The Fundamental Groupoid, John Wiley\& Sons, New York. Chichester. Brisbanc. Toronto.

\bibitem{dug} J. Dugundji, Topology, Universal Book Stall, New Delhi, (1990).

\bibitem{eki1} E. Ekici,  \newblock{On $e$-open sets, $DP^*$-sets and $DP\epsilon^*$-sets and decompositions of continuity}, \newblock{Arabian Journal for Science and Engineering}, { 33} (2A)(2008), 269-282.
\bibitem{eki2} E. Ekici,  \newblock{Some generalizations of almost contra-super-continuity}, \newblock{Filomat}, { 21} (2) (2007), 31-44.
\bibitem{eki3} E. Ekici,  \newblock{New forms of contra-continuity}, \newblock{Carpathian Journal of Mathematics}, {24} (1) (2008), 37-45.
\bibitem{eki4} E. Ekici, \newblock{On $ e^* $-open sets and $ (D,S)^* $-sets}, \newblock{Mathematica Moravica}, {13} (1) (2009), 29-36.
\bibitem{eki5} E. Ekici, \newblock{ A note on $ a $-open sets and $ e^* $-open sets}, \newblock{Filomat}, {22} (1) (2008), 89-96. 

\bibitem{mas} S. Mashhour,A. A. Allam, F. S. Mahmoud and F. H. Khedr, On Supra Topological spaces, Indian J. Pure and Appl.Math, $14(1983),502-510.$

\bibitem{mun} J. R. Munkres, Topology A First Course, Prentice-Hall, Inc. Englewood Clifs, New Jersey,(1975).

\bibitem{nja} O. Njastad, On Some Classes of Nearly Open Sets, Pacific Journal ofmathematics, Vol.15, No, 31965.

\end{thebibliography}
\end{document}